\documentclass[a4paper]{amsart}
\usepackage[ascii]{inputenc}
\usepackage[notcite]{showkeys}

\usepackage{amssymb}
\usepackage{mathrsfs}
\usepackage{hyperref}
\usepackage{cancel}

\usepackage{xcolor}

\usepackage[shortlabels]{enumitem}
\setlist[enumerate,1]{label={(\Alph*)}}
\setlist[enumerate,2]{label={(\alph*)}}
\setlist[enumerate,3]{label={$\bullet_{\arabic*}$}}

\newenvironment{PROOF}[2][\proofname.]
   {\begin{proof}[#1]\renewcommand{\qedsymbol}{\textsquare$_{\mathrm{#2}}$}}
   {\end{proof}}

\newtheorem{theorem}{Theorem}[section]

\newtheorem{claim}[theorem]{Claim}
\newtheorem{conclusion}[theorem]{Conclusion}

\newtheorem{observation}[theorem]{Observation}

\theoremstyle{definition}

\newtheorem{convention}[theorem]{Convention}

\newtheorem{definition}[theorem]{Definition}

\newtheorem{fact}[theorem]{Fact}
\newtheorem{hypothesis}[theorem]{Hypothesis}

\theoremstyle{remark}

\newtheorem{notation}[theorem]{Notation}

\newtheorem{remark}[theorem]{Remark}

\usepackage{graphicx}

\newcommand{\av}{\mathrm{av}}

\newcommand{\HJ}{\mathsf{HJ}}

\newcommand{\dom}{\mathrm{dom}}

\newcommand{\ob}{\mathrm{ob}}

\newcommand{\Sub}{\mathrm{Sub}}

\newcommand{\bfH}{\mathbf{H}}

\newcommand{\bfS}{\mathbf{S}}

\newcommand{\bfW}{\mathbf{W}}

\newcommand{\bfc}{\mathbf{c}}

\newcommand{\bfo}{\mathbf{o}}

\newcommand{\bfs}{\mathbf{s}}
\newcommand{\bft}{\mathbf{t}}
\newcommand{\bfu}{\mathbf{u}}

\newcommand{\bbN}{\mathbb{N}}

\newcommand{\bbR}{\mathbb{R}}

\newcommand{\mn}{\medskip\noindent}
\newcommand{\sn}{\smallskip\noindent}
\newcommand{\bn}{\bigskip\noindent}

\newcommand{\cE}{\mathscr{E}}

\newcommand{\clE}{\mathcal{E}}

\newcommand{\clU}{\mathcal{U}}

\newcommand{\gS}{\mathfrak{S}}

\newcommand{\gs}{\mathfrak{s}}
\newcommand{\gt}{\mathfrak{t}}
\newcommand{\gu}{\mathfrak{u}}

\newcommand{\eps}{\varepsilon}
\newcommand{\cl}{c\kern-.11ex \ell}
\newcommand{\lh}{{\ell\kern-.27ex g}}
\newcommand{\rest}{\restriction}

\newcommand{\LL}{\langle}
\newcommand{\RR}{\rangle}

\newcommand{\subref}[1]{$_{\mathrm{\texttt{=}}\mathsf{L{#1}}}$}

\usepackage{tcolorbox}

\newcommand{\overbar}[1]{\mkern 1.5mu\overline{\mkern-1.5mu#1\mkern-1.5mu}\mkern 1.5mu}

\newcommand*{\defeq}{\mathrel{\vcenter{\baselineskip0.5ex \lineskiplimit0pt\hbox{\scriptsize.}\hbox{\scriptsize.}}}=}

\newcommand{\leqdot}{\mathrel{\mathord{\leq}\!\!\raise
0.8 pt\hbox{$\cdot$}}}

\newcount\skewfactor
\def\mathunderaccent#1#2 {\let\theaccent#1\skewfactor#2
\mathpalette\putaccentunder}
\def\putaccentunder#1#2{\oalign{$#1#2$\crcr\hidewidth
\vbox to.2ex{\hbox{$#1\skew\skewfactor\theaccent{}$}\vss}\hidewidth}}

\newbox\noforkbox \newdimen\forklinewidth
\setbox0\hbox{$\textstyle\bigcup$}
\setbox1\hbox to \wd0{\hfil\vrule width \forklinewidth depth \dp0
                        height \ht0 \hfil}
\setbox\noforkbox\hbox{\box1\box0\relax}
\def\unionstick{\mathop{\copy\noforkbox}\limits}
\def\nonfork#1#2_#3{#1\unionstick_{\textstyle #3}#2}
\def\nonforkin#1#2_#3^#4{#1\unionstick_{\textstyle #3}^{\textstyle
    #4}#2}
\setbox0\hbox{$\textstyle\bigcup$}
\setbox1\hbox to \wd0{\hfil{\sl /\/}\hfil}
\setbox2\hbox to \wd0{\hfil\vrule height \ht0 depth \dp0 width
                                \forklinewidth\hfil}
\newbox\doesforkbox
\setbox\doesforkbox\hbox{\box1\box0\relax}
\def\nunionstick{\mathop{\copy\doesforkbox}\limits}

\def\fork#1#2_#3{#1\nunionstick_{\textstyle #3}#2}
\def\forkin#1#2_#3^#4{#1\nunionstick_{\textstyle #3}^{\textstyle
    #4}#2}

\newcommand{\stickT}{%
\setbox255=\hbox{\raise1ex\hbox{$\hspace{0.2pt}\,\bullet\,$}}
\mathord{\rlap{\hbox to\wd255{\hss\hbox{$|$}\hss}}
\box255}
}
\newcommand{\stickS}{%
\setbox255=\hbox{\raise0.6ex\hbox{$\scriptstyle\bullet$}}
\mathord{\rlap{\hbox to\wd255{\hss\hbox{$\scriptstyle|$}\hss}}
\box255}
}

\author[S. Shelah]{Saharon Shelah}
\address{Einstein Institute of Mathematics,
The Hebrew University of Jerusalem,
9190401, Jerusalem, Israel; and\\
Department of Mathematics,
Rutgers University,
Piscataway, NJ 08854-8019, USA}
\urladdr{https://shelah.logic.at/}
\thanks{First typed 2026-06-17.
The author would like to thank the Israel Science Foundation (ISF) for partially supporting
this research by
grant  2320/23,
and thanks Matt Grimes for the careful and beautiful typing.\\
References like (e.g.) [Sh:950, Th.0.2\subref{y5}] mean that \textsf{y5} is the internal label of Theorem 0.2 in the TeXfile of [Sh:950].
The reader should note that the version in my website is usually more up-to-date than the one in arXiv.
This is publication number
E117
in Saharon Shelah's list.
}

\makeatletter
\@namedef{subjclassname@2020}{\textup{2020} Mathematics Subject Classification}
\makeatother
\date{November 5, 2026}

\title{Hales-Jewett numbers revisited\\ E117}

\definecolor{altgreen}{rgb}{0,0.5,0}

\providecommand{\bfW}{\mathbf{W}}
\providecommand{\ob}{\mathrm{ob}}

\providecommand{\Sub}{\mathrm{Sub}}

\definecolor{whyblue}{rgb}{0,0.3,0.8}

\begin{document}
\begin{abstract}
Hales-Jewett theorem is well known.
We prove a primitive recursive bound, in fact an $\clE^5$ one, for the density version.

In preparation we deal with what was the main induction step in the old proof of Hales-Jewett as a result in its own right with better bounds. Also in improving our bounds.

\end{abstract}

\maketitle

\newpage
\centerline{\textbf{\underline{Annotated Content}}}

\mn
\textbf{\S0 \quad Introduction} \hfill p.\pageref{S0}

\bn
\textbf{\S0(B) \quad Preliminaries} \hfill p.\pageref{S0B}

\bn
\textbf{\S1 \quad The density version of Hales-Jewett, $\clE^5$-bound} \hfill p.\pageref{S1}

\newpage
\setcounter{section}{-1}
\section{Introduction}\label{S0}

Here we get a primitive recursive, even $\clE^5$, bound for $\HJ_\eps(\Lambda,r)$, Hales-Jewett with density, as a function of $\lfloor\frac1\eps\rfloor+|\Lambda|+r$.

The history of this problem is well-known; for the early parts, see the book [GraRotSpe90]. (For more recent history, see [She:679] or Wikipedia.)
Here in \S1 we get a bound in $\clE^5$.

The Density Hales--Jewett Theorem, stating that any dense subset of the grid $[\ell]^n$ contains a combinatorial line, was originally proven via ergodic theory by Furstenberg and Katznelson [FurKat91]. A landmark quantitative elementary combinatorial proof was  discovered (considerably later) by Polymath 1 [Pol12], a massively collaborative project.
Their proof gives bounds in terms of a version of the Ackermann function.
Here we give a primitive recursive bound;
in fact in $ {\mathscr E} ^5$; it is two levels above exponentiation in the Grzegorczyk hierarchy.

Recent work in [GolMirShe:1182], among other related results, connects van der Waerden's Theorem and Hales-Jewett.

\bn
 Nowadays I will surely be asked how much AI helped me:
\begin{enumerate}
    \item [(a)] in the \LaTeX\ typing, it helped translating my handwriting to \LaTeX;
    \item [(b)] it helped in proofreading (pointing out errors, inconsistencies etc.);
    \item [(c)] we shall discuss \S3 when making it public.
\end{enumerate}

Recall the Hales-Jewett Theorem [HalJew63].
\begin{theorem}[\textbf{Hales--Jewett Theorem}]\label{z0}
For all natural numbers $n,p \ge 2$ there exists an integer $\HJ(n,p)$ such that for all {linear orders} $M$ of cardinality $\geq \HJ(n,p)$, $\Lambda$ of cardinality
$n$, $C$ of cardinality $p$, and $\bfc : {}^M\!\Lambda \to C$, there exists a combinatorial line (see below) on which $\bfc$ is constant.

(Think of $\bfc$ as a coloring: $\Lambda$ is our alphabet, and $C$ the set of colors.)
\end{theorem}

\bn
\begin{definition}\label{z1}
1) A \emph{combinatorial line} in ${}^M\!\Lambda$ is a set of the form
$$
\{ \rho \in {}^M\!\Lambda : \rho \supseteq \eta \text{ and $\rho \rest u$ is constant}\},
$$
for some non-empty $u \subseteq M$ and $\eta \in {}^{M\setminus u}\Lambda$.

\mn

2) Generalizing this, we say $\Gamma$ is a \emph{subspace} of ${}^M\!\Lambda$ of dimension $\ell$ (or $\ell$-subspace) \underline{when} it is of the form
$$
\big\{ \rho \in {}^M\!\Lambda : \rho \supseteq \eta, \text{ and $m < \ell \Rightarrow \rho \rest u_m$ is constant} \big\}
$$
for some sequence $\LL u_m : m < \ell\RR$ of pairwise disjoint non-empty subsets of $M$ and $\eta : M \setminus \bigcup\limits_{m<\ell} u_m \to \Lambda$.

\mn
 3) Let $\HJ(n,r,p)$ be defined as in \ref{z0}, but for a monochromatic $r$-subspace (so $\HJ(n,1,p) \defeq \HJ(n,p)$).
\end{definition}

\bigskip
The following statement reduces to the van der Waerden Theorem for $m=1$. For $m > 1$
this is the Gallai-Witt Theorem
$=$ the multidimensional van der Waerden Theorem.
\begin{theorem}\label{z5}
Let $n \ge 1$ and $m,p \ge 2$. Then there exists a number $k$ satisfying the following:
\begin{enumerate}
    \item [$(*)_{n,p}^{k,m}$] If $\LL A_\ell : \ell < p\RR$ is a partition of ${}^nk$, \underline{then} there exist $\bar a \in {}^nk$, $b_0,\ldots b_{n-1} \in k \setminus \{0\}$, and $\ell^* < p$ such that
    $$
    \big\{\bar a + \big\LL \eta(i)b_i : i < n\big\RR : \eta \in {}^nm \big\} \subseteq A_{\ell^*}.
    $$
\end{enumerate}

(This can be deduced from $\HJ(n,k,p)$; see [GraRotSpe90]. The proof applies to the density version as well,
so our bounds follow.)
\end{theorem}

\bn
\begin{definition}\label{z8}
Assume $\eps \in (0,1]$. Let $\HJ_\eps(\Lambda,r)$ be the minimal natural number $k$ satisfying the following implication.
\begin{quotation}
    If $X \subseteq {}^k\!\Lambda$ has density $\ge \eps$ (i.e.\ $|X|\,/\,|{}^k\!\Lambda| \ge \eps$) \underline{then} $X$ contains an $r$-dimensional subspace of ${}^k\!\Lambda$.
\end{quotation}
\end{definition}

\bn
\begin{fact}\label{z11}
If $k = \HJ_\eps(n,r)$ and $\eps\le\frac1p$ then $k\ge\HJ(n,r,p)$.
\end{fact}

\begin{definition}\label{y14}
$\clE^n$ denotes the $n$-th class of the Grzegorczyk hierarchy of primitive recursive functions. So a function is in $\clE^3$ (the elementary functions) iff it is bounded by a fixed number of iterated exponentiations, it is in $\clE^4$ iff it is bounded by an iteration of exponentiation (a tower of exponents of height a fixed polynomial, or more generally an $\clE^3$ function, of the argument), and $\clE^5$ allows in addition an iteration of the tower function.
\end{definition}

\bn
\begin{notation}\label{y11}
1) Here $i,j,k,\ell,m,n,p,r,s,t$ will denote natural numbers; equivalently, finite ordinals. (So $n$ is identified with $\{i : i < n\}$.)

\mn
2) $\Lambda$ is a finite alphabet (mostly fixed). $\alpha,\beta,\gamma$ will denote members of $\Lambda$, usually; we fix $\alpha$ and $\beta$.

\mn
3) $\eps$, $\zeta$, $\xi$ and $\varkappa$, $w$, $W$ will be real numbers in the interval $(0,1]_\bbR$.

\mn
4) $\bfs$, $\bft$, and $\bfu$ will denote so called systems (see Definition \ref{d8}). $\gs$, $\gt$, $\gu$ will denote pre-systems, and $\bfo\in\ob(\gs)$.

\mn
5) $\exp_m(n) = \exp(m,n) \defeq m^n$ (We will use this to avoid deeply nested superscripts.)

\mn
5A) To this end, we let $\exp_m^0\!(n) \defeq n$ and $\exp_m^{k+1}(n) \defeq \exp_m(\exp_m^k\!(n))$.

\mn
6) ${}^A\!B$ is the set of functions from $A$ to $B$.

\mn
7) $m^+$ denotes the successor of $m$ (so it is shorthand for $m+1$, as these are natural numbers).

\mn
8) We shall use $\{k,\ell\}$ with $k < \ell$ to denote a member of $[m]^2 \defeq \{A \subseteq m : |A| = 2\}$. We will endeavor to make sure that when these two letters are used together it is only in this context.

\mn
9) If we write (e.g.) $\HJ(\Lambda,C)$, we of course mean $\HJ\big(|\Lambda|,|C|\big)$, and if in the definition there is no natural number satisfying the demand then it means $\infty$.

\mn
 10) $N$, $\mathcal M$, $\mathcal N$ and $M$ are finite linear orders and for $a \in M \cup \{\infty\}$ let $M_{<a} = \{b \in M : b <_M a\}$.

\mn
 11) Let $M = M_1 \,\dot\cup\, M_2$ mean: $M_1\cap M_2 = \emptyset$, $(a\in M \Leftrightarrow a\in M_1 \vee a\in M_2)$ and ${<_M}\rest M_\iota = {<_{M_\iota}}$ for $\iota=1,2$.
\end{notation}

\subsection{Explanation}\label{S0A}
Here, on the one hand, we give new bounds, and on the other hand we give some bounds to some relatives; usually they are inspired by the main inductive step in the proof of [She:329, \S1].

We shall first try to explain the earlier proofs, and only then the present one.

\mn
 \textbf{First, the earlier proofs.}

The proof of the Hales-Jewett Theorem [HalJew63] uses double induction: first on $n$ (the cardinality of the alphabet $\Lambda$), and then on $p$ (the number of colors). Eventually they get an upper bound in terms of the Ackermann function, which is not primitive recursive.

The base idea of [She:329] is that we fix the number of colors $|C|$ and perform our induction only on $n \defeq |\Lambda|$.

Let us try to describe the inductive step in the proof (that is, from $n$ to $n+1$ letters). So we let $r \defeq \HJ(n,p)$ and fix $\Lambda$ of cardinality $n+1$, $C$ of cardinality $p$, and $\alpha \neq \beta \in \Lambda$. We look for a finite linear order $M$ which is large enough for our purposes (i.e.\ $|M| \gg r$) and consider a coloring $\mathbf{c} : {}^M\Lambda \to C$.

Let $\overbar N = \LL N_i : i < r\RR$ be a partition of $M$ into intervals; the intention is to use each $N_i$ as a place for one letter. That is, we will try to find $\xi = \big\langle (k_i,\ell_i,N_i) : i < r\big\rangle$ with $k_i < \ell_i$ both from\footnote{Actually, we allow $\ell_i$ to be above all elements of $N_i$.} $N_i$, and restricting ourselves to functions $\eta \in {}^M\Lambda$ of the form
$$
i < r \wedge j \in N_i \ \Rightarrow\ \eta(j) \defeq
\begin{cases}
\alpha & \text{if } j < k_i\\
\gamma_i & \text{if } j \in [k_i,\ell_i)\\
\beta & \text{if } j \ge \ell_i.
\end{cases}
$$
This was done such that the coloring does not change if we replace some of the $\alpha$-s by $\beta$-s (or vice versa): we reduce this case of the problem to one with one letter less, hence can use the induction hypothesis, that is, the choice of $r$ gives the desired conclusion.

So $k_i$ and $\ell_i$ were chosen by downward induction on $i < r$, and naturally we required $|N_i| > |C|^{\sum_{j<i}|N_j|}$, so for our lemma we had a tetration function as a bound. Hence for $\HJ$ we had gotten a bound in $\clE^5$.

\textbf{Second, our present aim.}

Here we have three aims concerning $\HJ(\Lambda,C)$ (not disjoint):
\begin{enumerate}
    \item [(A)] improve the bound,
    \item [(B)] deal with the density versions,
    \item [(C)] get good bounds for relatives, mainly for the inductive step of [She:329].
\end{enumerate}
Sometimes the bounds we get for a claim are superseded if we consider only applying them to $\HJ_\eps(\Lambda,r)$, but not so if we consider the proof and the specific results it gives on the way.

\mn
 \textbf{Third, the present work:} We first concentrate on the major inductive step. We have $\gs$ consisting of $\mathcal M$ and $M = \sum\{M_a : a\in\mathcal M\}$, and the combinatorial space is $V_\gs = {}^M\Lambda$. The density enters by expanding $\gs$ to $\bfs$ by finite $A_\bfs\ne\emptyset$ and $A^\bfs_\eta\subseteq A_\bfs$ for $\eta\in V_\gs$. It is natural to consider $\eps$ = the average of $|A^\bfs_\eta|/|A_\bfs|$, but (wlog, \ref{d65}) we work with $\eps = \min(\bfs)$ = the minimal $|A^\bfs_\eta|/|A_\bfs|$.

We like to find a subspace of dimension $|\mathcal M|$ which is colour blind for fixed $\alpha\ne\beta\in\Lambda$.

As in [She:329] the subspace we consider is represented by $\LL(k_a,\ell_a) : a\in\mathcal M\RR$ where $k_a<\ell_a<|M_a|$ (really $\ell_a\le|M_a|$ suffices). It consists of considering $\eta\in V_\gs$ such that for every $a\in\mathcal M$, $\eta\rest M_a$ is $\alpha$ on the first $k_a$ elements, $\beta$ for the last $|M_a|-\ell_a$ elements and constant on the rest.

It had seemed most natural to decrease the density $\eps$ during the induction, but eventually a different road is taken: we fix $\eps$ but for $t_*$ ask about being so called $(\eps,t)$-good, related to having a subspace of dimension $t\le t_*$, intersecting the $A_\eta$ where in the free subsets (of $M$) we put $\alpha$ or we put $\beta$. This is carried in the induction on $|\Lambda|$.

We forgot to say, the pigeon hole principle is replaced by Corr\'adi lemma (\ref{x4}).

\subsection{Preliminaries}\label{S0B}

We recall some needed facts.

\bn
\begin{notation}\label{x1}
1) For $A$ a finite non-empty set, let $\mu_A(-)$ be the $A$-\emph{density} function. I.e.\ $\mu_A(B) \defeq |B| \,/\, |A|$ for $B \subseteq A$.

(We may write $\mu$ instead of $\mu_A$ when $A$ is clear from the context.)

\mn
2) Let $A^{[\ge\eps]} \defeq \{B \subseteq A : \mu_A(B) \ge \eps\}$.

\end{notation}

\bn
 Recalling the Corr\'adi lemma:
\begin{fact}\label{x4}
1) We have `(A) $\Rightarrow$ (B),' where
\begin{enumerate}
    \item
    \begin{enumerate}
        \item $n \geq 2$ and $\eps \in \big[\frac1n ,1\big]_\bbR$.
\sn
        \item $A$ is a finite non-empty set and $\mu= \mu_A$ is as in \ref{x1}.
\sn
        \item For $i<n$ we have $A_i \subseteq A$ with $\mu(A_i) \geq \eps$.
\sn
        \item $0 < \zeta \leq \dfrac{\eps(\eps n - 1)}{n-1}$
    \end{enumerate}
\sn
    \item There exist $k < \ell < n$ such that $\mu(A_k \cap A_\ell) \geq \zeta$.
\end{enumerate}

\mn
2) In \ref{x4}(1) we can weaken  clause (c) to
$ \dfrac{1}{n} \Sigma _{k < n } \mu (A_k)
   \ge \eps$
\end{fact}

\sn
\begin{PROOF}{\ref{x4}}
1) This is the Corr\'adi Lemma [Cor69a], [Cor69b].

\mn
2) A variant of Corr\'adi lemma, see [Juk11, Chapter 2, Lemma 2.2].
\end{PROOF}

\bn
\begin{observation}\label{x7}
1) If $n \ge 1 + \lceil\frac{1}{\eps}\rceil$ and non-zero $\zeta \le \eps^3$ then $n$, $\eps$ and $\zeta$ will satisfy the conditions of \emph{\ref{x4}(1)(A)(a),(d)}.

\mn
 2) If $k \ge 1$, $n \ge 1 + \left\lceil\dfrac{1}{\eps^{k+1}}\right\rceil$ and $0 < \zeta \le \eps^2\big(1 - \eps^{k}\big)$, then the same conclusion holds.
\end{observation}

\sn
\begin{PROOF}{\ref{x7}}
1) Obviously clause \ref{x4}(1)(A)(a) holds. For clause (A)(d), recall that
$$
\frac{\eps(\eps n - 1)}{n-1} = \eps \cdot \frac{\eps (n-1) - (1-\eps)}{n-1} = \eps^2 - \frac{(1-\eps)\eps}{n-1}.
$$
The last expression is increasing in $n$ (even for real-valued $n$), hence it is
$$
\ge \eps^2 - \frac{(1-\eps)\eps}{(1+\frac{1}{\eps})-1} = \eps^2 - (1-\eps)\eps^2 = \eps^3,
$$
so `$\zeta \le \eps^3$' does indeed suffice.

\mn
 2) Similarly.
\end{PROOF}

\newpage
\section{An $\clE^{5}$ bound for $\HJ(\Lambda,r)$\protect\footnote{In various places we could get better bounds, but immaterial here.}}\label{S1}

\bn
\begin{convention}\label{d0}
1) Recall $\gs,\gt,\gu$ vary on $\gS$ and $\bfs,\bft,\bfu$ vary on $\bfS$.

2) Superscripts/subscripts like $\gs$, $\bfs$ and even $\eps$ may be omitted when clear from the context.

3) Always $\Omega\subseteq{}^M\Lambda$ and $\Upsilon\subseteq{}^{\mathcal M}\Lambda$.
\end{convention}

\bn
\begin{hypothesis}\label{d5}
1) $\Lambda$ a finite alphabet.

2) $\alpha\ne\beta\in\Lambda$, fixed till \ref{d70}.
\end{hypothesis}

\bn
\begin{definition}\label{d8}
1) Let $\gS = \gS_{\alpha,\beta}$ be the set of pre-systems $\gs$ meaning (so $M = M_\gs$, $M_a = M^\gs_a$, etc.):
\begin{enumerate}
    \item [(a)] $\gs = (M,\mathcal M,\bar M,\clU)$,
    \item [(b)] $\mathcal M$ a finite linear order,
    \item [(c)] $\bar M = \LL M_a : a\in\mathcal M\RR$ and $M = \sum\{M_a : a\in\mathcal M\}$,
    \item [(d)] $\clU = \{a : |M_a| = 1\}$,
    \item [(e)] $a\in\mathcal M\setminus\clU\Rightarrow|M_a|\ge2$.
\end{enumerate}

1A) For $\gs\in\gS$:

\quad(a) let $V_\gs = {}^M\Lambda$;

\quad(b) for $\mathcal N\subseteq\mathcal M$ let $\gs\rest\mathcal N = \big(\bigcup_{a\in\mathcal N}M_{\gs,a},\ \bar M\rest\mathcal N,\ \clU_\gs\cap\mathcal N\big)$;

\quad(c) $\Upsilon_\gs = {}^{\mathcal M}\Lambda$;

\quad(d) $\delta = \delta_{\gs,\eps} : \mathcal M\to[0,1]_\bbR$ is defined by: if $a\in\mathcal M$ and $r_\gs(a) := |\{b\in\mathcal M_\gs : b<a\}|$ then $\delta_\gs(a) = \exp_\eps\exp_3(r_\gs(a))$.

2) $\ob(\gs)$ is the set of $\bfo = (\mathcal N,\cE,\rho) = (\mathcal N_\bfo,\cE_\bfo,\rho_\bfo)$ such that $\mathcal N_\bfo\subseteq\mathcal M$, $\cE_\bfo$ an equivalence relation on $\mathcal N_\bfo$, $\rho_\bfo\in{}^{\mathcal M\setminus\mathcal N}\Lambda$ and $a\in\mathcal N_\bfo\cap\clU_\gs\Rightarrow a/\cE_\bfo = \{a\}$.

2A) Let $\mathrm{sob}(\gs)$ be the set of $\bfo\in\ob(\gs)$ such that $\cE_\bfo$ is the equality.

2B) For $\bfo\in\ob(\gs)$ and $\mathcal M_*\subseteq\mathcal M_\gs$ let $\bfo\rest\mathcal M_* = \big(\mathcal N_\bfo\cap\mathcal M_*,\ \cE_\bfo\rest\mathcal M_*,\ \rho_\bfo\rest(\dom(\rho_\bfo)\cap\mathcal M_*)\big)$.

3) For $\bfo\in\ob(\gs)$ let
\begin{itemize}
    \item [$\bullet_1$] $\Upsilon_\bfo = \{\nu\in{}^{\mathcal M}\Lambda : \nu$ extends $\rho_\bfo$ and for each equivalence class $u$ of $\cE_\bfo$, $\nu\rest u$ is constantly $\alpha$ or constantly $\beta\}$,
    \item [$\bullet_2$] $\dim(\bfo) = |\mathcal N_\bfo/\cE_\bfo|$ ($\dim$ is for dimension),
    \item [$\bullet_3$] $\mathrm{bas}(\bfo) = |\mathcal N_\bfo|$ ($\mathrm{bas}$ is for basis),
    \item [$\bullet_4$] $\Upsilon^+_\bfo = \{\nu\in{}^{\mathcal M}\Lambda : \nu$ extends $\rho_\bfo$ and is constant on each $\cE_\bfo$-equivalence class$\}$,
    \item [$\bullet_5$] $\ob_t(\gs) = \{\bfo\in\ob(\gs) : \dim(\bfo)\le t\}$.
\end{itemize}

4) Let $\bfH = \bfH_\gs$ be the set of $\mathbf h\in\prod\{[|M_a|+1]^2 : a\in\mathcal V_\gs\}$ where $\mathcal V_\gs = \mathcal M_\gs\setminus\clU_\gs$, and we may write $\mathbf h$ as
$$\LL(k_a,\ell_a) : a\in\mathcal V_\gs\RR = \LL(k_{\mathbf h,a},\ell_{\mathbf h,a}) : a\in\mathcal V_\gs\RR$$
where $k_a<\ell_a\le|M_a|$.

5) For $\mathbf h\in\bfH_\gs$, $\bfo\in\ob(\gs)$, $\rho\in\Upsilon_\bfo$ let $\eta = \eta[\rho,\mathbf h] = \eta_\gs[\rho,\mathbf h]$ be the $\eta\in V_\gs$ such that
if $a\in \mathcal M$ then $\eta\rest M_a$ is $\varrho_{k_{\mathbf h,a},\ell_{\mathbf h,a},\rho(a)}[M_a]$, see \ref{d11} below.

6) Let $\Omega_{\gs,\mathbf h,\bfo} = \{\eta_\gs[\rho,\mathbf h] : \rho\in\Upsilon_\bfo\}$ and $\Omega^+_{\gs,\mathbf h,\bfo} = \{\eta_\gs[\rho,\mathbf h] : \rho\in\Upsilon^+_\bfo\}$.

Let $\Omega_{\gs,\mathbf h} = \bigcup\{\Omega_{\gs,\mathbf h,\bfo} : \bfo\in\ob(\gs)\}$, and $\bfW_\gs$ $= \bigcup\{\Omega_{\gs,\mathbf h} : \mathbf h\in\bfH_\gs\}$.

7) Let $\bfW_\bfo = \bfW_{\gs,\bfo} = \bigcup\{\Omega_{\gs,\mathbf h,\bfo} : \mathbf h\in\bfH_\gs\}$.
\end{definition}

\bn
\bn
\begin{definition}\label{d11}
1) Assume $M$ is a linear order
\begin{enumerate}
    \item [(a)] if $k \le |M|$ then $\varrho_k = \varrho_k[M]$ is the member of ${}^M\Lambda$ defined by (for $a \in M$)
    \begin{itemize}
        \item $\varrho_k(a) = \alpha$ \underline{if} $|\{b : b <_M a\}| < k$,
        \item $\varrho_k(a) = \beta$ \underline{if} otherwise,
    \end{itemize}
    \item [(b)] if $k < \ell \le |M|$ and $\gamma \in \Lambda$ then $\varrho_{k,\ell,\gamma} = \varrho_{k,\ell,\gamma}[M]$ is defined by (for $a \in M$)
    \begin{itemize}
        \item $\varrho_{k,\ell,\gamma}(a) = \alpha$ \underline{if} $|\{b : b <_M a\}| < k$,
        \item $\varrho_{k,\ell,\gamma}(a) = \beta$ \underline{if} $|\{b : b <_M a\}| \ge \ell$,
        \item $\varrho_{k,\ell,\gamma}(a) = \gamma$ \underline{if} otherwise.
        \item note that $\varrho_{k,\ell,\alpha} = \varrho_\ell$ and $\varrho_{k,\ell,\beta} = \varrho_k$.
    \end{itemize}
    \item [(c)] Above we may replace $M$ by $\bfs$, meaning $M_\bfs$, similarly for $\gs$.
\end{enumerate}
\mn
2) If $\alpha,\beta$ are not fixed we may write $\varrho_k[M,\alpha,\beta]$, $\varrho_{k,\ell,\gamma}[M,\alpha,\beta]$; if also $\Lambda$ is not fixed we may write $\varrho_k[M,\alpha,\beta,\Lambda]$, $\varrho_{k,\ell,\gamma}[M,\alpha,\beta,\Lambda]$.

\mn
 2A) Let $\Omega^+_M = \Omega^+[M] = \{\varrho_{k,\ell,\gamma}[M] : k<\ell\le|M| \text{ and } \gamma\in\Lambda\}$.

3) Assume $\gt = \gs\rest(\mathcal M_\gs\setminus\{a\})$. For $k\le|M_{\gs,a}|$ and $\Omega \subseteq V_\gt$ let $\Omega+\varrho_k = \{\eta\cup\varrho_k[M_{\gs,a}] : \eta\in\Omega\}$.
\end{definition}

\bn
\begin{definition}\label{d14}
1) We say $\bfs$ is a system or $\bfs\in\bfS = \bfS_{\alpha,\beta}$ when $\bfs$ consists of:
\begin{enumerate}
    \item [(a)] $\gs = \mathrm{pre}(\bfs)$ a pre-system,
    \item [(b)] $A_\bfs$ a finite non-empty set,
    \item [(c)] $\bar A = \LL A_\eta : \eta\in V_\gs\RR$ where $A_\eta\subseteq A_\bfs$.
\end{enumerate}

2) If $\Omega\subseteq V_\gs$ then we let $A^\bfs_\Omega = \bigcap\{A^\bfs_\eta : \eta\in\Omega\}$.
\end{definition}

\bn
\begin{definition}\label{d17}
For $\bfs\in\bfS$:

1) Let $\av(\bfs) = \sum\{\mu_\bfs(A^\bfs_\eta) : \eta\in V_\bfs\}/|V_\bfs|$.

2) Let $\min(\bfs) = \min\{\mu_\bfs(A^\bfs_\eta) : \eta\in V_\bfs\}$.

3) We say $\bfs$ is $\eps$-good when: if $\bfo\in\ob(\bfs)$, $\mathcal N_\bfo\subseteq\clU_\bfs$ and $\mathbf h\in\bfH_\bfs$ then $\mu_\bfs(A^\bfs_{\Omega_{\gs,\mathbf h,\bfo}})\ge\exp_\eps(\exp_3(\dim(\bfo)))$.

4) We say $\bfs$ is $(\eps,t)$-good when: above we restrict ourselves to $\bfo$ such that $\dim(\bfo)\le t$.
\end{definition}

\bn
\begin{remark}\label{d20}
1) We may restrict ourselves to $\Omega_\gs\subseteq V_\gs$; this is natural in the induction on $|\Lambda|$. The saving seems insignificant.
\end{remark}

\bn
\begin{observation}\label{d23}
Note that if $\clU_\gs = \emptyset$ and $\min(\bfs)\ge\eps$ then $\bfs$ is $\eps$-good.
\end{observation}

\bn
\bn
\begin{definition}\label{d35}
1) For $\bfs$ a system and $\bfo\in \mathrm{sob}(\mathrm{pre}(\bfs))$ such that $\mathcal N_\bfo\cap\clU_\gs = \emptyset$, and $\mathbf h\in\prod\{[M_a]^2 : a\in\mathcal N_\bfo\}$, let $\Sub_{\mathbf h,\bfo}(\bfs)$ be the set of systems $\bft$ satisfying the following.
\begin{enumerate}
    \item [(A)] $\mathcal M_\bft = \mathcal M_\gs$ and $M_{\bft,a}$ is $M_{\gs,a}$ if $a\in\mathcal M_\gs\setminus\mathcal N_\bfo$, and is $\{M_{\gs,a}\}$ otherwise.
    \item [(B)] $\bullet$ $\clU_\bft = \clU_\gs\cup\mathcal N_\bfo$,

    $\bullet$ we define a function $h = h_{\mathbf h,\bfo} = h_{\gs,\mathbf h,\bfo}$ such that $h : V_\bft\to V_\bfs$ is defined by: $h(\nu)\rest M_{\gs,a}$ is
    \begin{itemize}
        \item $\nu\rest M_{\gs,a}$ \underline{if} $a\in\mathcal M_\gs\setminus\mathcal N_\bfo$,
        \item $\varrho_{k,\ell,\gamma}$ \underline{if} $a\in\mathcal N_\bfo$, $\mathbf h(a) = (k,\ell)$ and $\nu(M_{\gs,a}) = \gamma$
    \end{itemize}
    (recalling $M_{\bft,a}$ is the singleton $\{M_{\gs,a}\}$).
    \item [(C)] $A^\bft_\eta\subseteq A^\bfs_{h(\eta)}$ for all $\eta\in V_\bft$.
    \item [(D)] $A_\bft = A_\bfs$.
\end{enumerate}

2) Let $\Sub_\bfo(\bfs) = \bigcup\{\Sub_{\mathbf h,\bfo}(\bfs) : \mathbf h \text{ as above}\}$ and $\Sub(\bfs) := \bigcup\{\Sub_\bfo(\bfs) : \bfo\in\ob(\bfs)\}$.

3) $\Sub_\bfo(\gs)$ and $\Sub(\gs)$ are defined similarly, omitting clauses (1)(C),(D).
\end{definition}

\bn
\begin{definition}\label{d36}
1) For $\gs,\gt\in\gS_{\alpha,\beta}$ let $\gs\to_{\min}(\gt)_{(\eps ,t)}$ mean
\begin{enumerate}
    \item [(a)] $\gs,\gt$ are pre-systems,
    \item [(b)] $0<\eps\le1$,
    \item [(c)] if $\bfs$ is an $(\eps ,t)$-good system expanding $\gs$ with $\min(\bfs)\ge\eps$ \underline{then} there exists a system $\bft$ with $\mathrm{pre}(\bft) = \gt$ and $\min(\bft)\ge\eps$ which is isomorphic to an $(\eps ,t)$-good member $\bfu$ of $\Sub(\bfs)$.
\end{enumerate}

2) For $\gs,\gt\in\gS$ let $\gs\to_{\av}(\gt)_{\eps _*,\eps ,t}$ be defined as in \ref{d36}(1) replacing $\min$ by $\av$ in clause (c), twice.

3) For $\gs,\gt\in\gS$ let $\gs\to_{\rm avm}(\gt)_{\eps_*,\eps,t}$ be defined as in \ref{d36}(1) replacing $\min(\bfs)\ge\eps$ by $\av(\bfs)\ge\eps_*$ in clause (c).
\end{definition}

\bn
\begin{claim}\label{d40}
We have $\gs\to_{\min}(\gt)_{\eps,t}$ when
\begin{enumerate}
    \item [$\boxplus$] (a) $\mathcal M_\gs = \mathcal M = \mathcal M_\gt = \clU_\gt$, and $t\le$ $r = |\mathcal M|$,

    (b) $0<\eps <1$,

    (c) $|M_{\gs,a}| = m>1$ so $\clU_\gs = \emptyset$,

    (d) $s = \big\lceil\frac1{\exp_\eps\exp_3(t)}\big\rceil+1$,

    (e) $m\to(sr)^{2r}_c$ where $c = 2^{|\ob(\gs)|}$, so $c\le\exp_2\big((2(t+|\Lambda|))^{|\mathcal M|}\big)$.
\end{enumerate}
\end{claim}

\bn

\begin{remark}\label{d42}
In \ref{d40}(e) we could use $m+1$.
\end{remark}

\begin{PROOF}{\ref{d40}}
By Definition \ref{d36}, let
\begin{enumerate}
    \item [$(*)_0$] $\bfs$ be an expansion of $\gs$ such that $\min(\bfs)\ge\eps$, and we should find suitable $\bfu$ and $\bft$.
\end{enumerate}
We define a function $F$ by
\begin{enumerate}
    \item [$(*)_1$] (a) $\dom(F) = \bfH_\gs$,

    (b) $F(\mathbf h)$ is the set of $\bfo\in\ob _t(\gs)$ such that $\mu_\bfs(A^\bfs_{\Omega_\bfo})\ge\exp_\eps(\exp_3(\dim(\bfo)))$,
    \item [$(*)_2$] (a) $\mathrm{Rang}(F)$ has $\le|\mathscr P(\ob _t(\gs))|$ elements,

    (b) this number is $\le \exp_2((2(t+|\Lambda|))^{|\mathcal M|})$,

    [Why? For $\bfo\in\ob_t(\gs)$ let $f_\bfo$ be the function with domain $\mathcal M_\gs$ such that
    \begin{itemize}
        \item if $a\in\mathcal N_\bfo$ then $f_\bfo(a) = (1,|\{b/\cE_\bfo : \min(b/\cE_\bfo)<\min(a/\cE_\bfo)\}|)$,
        \item if $a\in\mathcal M_\gs\setminus\mathcal N_\bfo$ then $f_\bfo(a) = (2,\rho_\bfo(a))$.
    \end{itemize}
    So the mapping $\bfo\mapsto f_\bfo$ is one-to-one and the number of such $f$-s is $\le\sum\{(s+|\Lambda|)^{|\mathcal M|} : s\le t\}\le (2(t+|\Lambda|))^{|\mathcal M|}$.]
    \item [$(*)_3$] there is a set $S\subseteq m$ such that

    (a) $S\subseteq m$,

    (b) $|S| = sr$,

    (c) $F\rest\bfH_S$ is constant where

    $\bfH_S = \{\mathbf h\in\bfH_\gs : \mathrm{rang}(\mathbf h)\subseteq S$ and for every $a<b$ from $\mathcal M$ we have $\max(\mathbf h(a))<\min(\mathbf h(b))\}$.

    [Why? By Ramsey theorem, i.e.\ by clause (e) of the claim.]
    \item [$(*)_4$] Choose $\bar S$ such that

    (a) $\bar S = \LL S_a : a\in\mathcal M\RR$,

    (b) the $S_a$-s are pairwise disjoint with union $S$,

    (c) $|S_a| = s$ and $a<_{\mathcal M}b\Rightarrow\max(S_a)<\min(S_b)$,
    \item [$(*)_5$] (a) for $a\in\mathcal M$ let $k_a = \min(S_a)$, $\ell_a = \max(S_a)$,

    (b) let $\mathbf h = \LL(k_a,\ell_a) : a\in\mathcal M\RR$,

    (c) we define $\bft$ expanding $\gt$ by: $A_\bft = A_\bfs$ and if $\rho\in V_\gt$ then $A^\bft_{\rho} = A^\bfs_{\eta[\rho,\mathbf h]}$.
\end{enumerate}
The main point to check is
\begin{enumerate}
    \item [$(*)_{6}$] $\bft$ is $(\eps ,t)$-good.
\end{enumerate}
This is a condition for every $\bfo\in\ob(\gt)$ with $\dim(\bfo)\le t$, and we do it by induction on $\dim(\bfo)$ using \ref{x4}(1).

First, $\dim(\bfo) = 0$: then the condition is $\mu_\bft(A^\bft_{\rho_\bfo})\ge\eps$, which holds as $A^\bft_{\rho_\bfo} = A^\bfs_{\eta[\rho_\bfo,\mathbf h]}$ and $\min(\bfs)\ge\eps$, which holds by $(*)_0$.

Second, in details, assume $\bfo$ with $\dim(\bfo) = n+1\le t$, let $\mathcal N$ be one of the $\cE_\bfo$-equivalence classes, for $a\in\mathcal N$ and $i<s$ let $k_{a,i}$ be the $i$-th member of $S_a$ and $\varrho_{a,i} = \varrho_{k_{a,i}}[M_a]$, and $\varrho_i = \bigcup\{\varrho_{a,i} : a\in\mathcal N\}$. Next we define $\bfo_i\in\ob(\gs)$ as follows:
\begin{itemize}
    \item $\mathcal N_{\bfo_i} = \mathcal N_\bfo\setminus\mathcal N$,
    \item $\cE_{\bfo_i} = \cE_\bfo\rest\mathcal N_{\bfo_i}$,
    \item $\rho_{\bfo_i} = \rho_\bfo\cup\varrho_i$.
\end{itemize}
By the induction hypothesis
\begin{enumerate}
    \item [$\circledast$] $\mu(A^\bfs_{\Omega_{\bfo_i,\mathbf h}})\ge\exp_\eps\exp_3(\dim(\bfo_i)) = \exp_\eps (\exp_3(n))$.
\end{enumerate}
As $n<t$, by the choice of $s$ in clause (d) of $\boxplus$, by \ref{x4}(1) for some $i<j<s$, $\mu(A^\bfs_{\Omega_{\bfo_i,\mathbf h}}\cap A^\bfs_{\Omega_{\bfo_j,\mathbf h}})\ge(\exp_\eps\exp_3(n))^3 = \exp_\eps\exp_3(n+1)$. By the choice of $S$ this holds for every $i<j<s$. This suffices.
\end{PROOF}

\bn

\bn
\begin{definition}\label{d62}
For $\gs\in\gS$ and $u\subseteq\bbN$ finite $\ne\emptyset$ (e.g.\ $k$) let $\gt = u\times\gs$ be the following member of $\gS$:
\begin{enumerate}
    \item [(a)] $\mathcal M_\gt = \mathcal M_\gs$ and $\clU_\gt = \clU_\gs$,
    \item [(b)] $M^\gt_a = u\times M^\gs_a$ for $a\in\mathcal M_\gs\setminus\clU_\gs$, ordered lexicographically,
    \item [(c)] $M^\gt_a = M^\gs_a$ for $a\in\clU_\gs$.
\end{enumerate}
\end{definition}

\bn
 In \ref{d65} no harm in assuming $\clU_\gs = \emptyset$.

\begin{claim}\label{d65}
If (A) then (B) $+$ (C), where:
\begin{enumerate}
    \item [(A)] (a) $\gs\in\gS$,

    (b) $0<\eps<\eps_*<1$,

    (c) $\zeta = \dfrac{\eps_*-\eps}{|V_\gs|-1}$, $k = \lfloor\frac1\zeta \rfloor+1$, and $\gu = k\times\gs$,

    (d) $\bfu\in\bfS$ expands $\gu$ and $\av(\bfu)\ge\eps_*$;
    \item [(B)] there are $\ell$, $\varrho_*$, $\bft = \bft_\ell$ such that
    \begin{enumerate}
        \item [(a)] $\ell<k$,
        \item [(b)] $\varrho_*\in V_{\ell\times\gs}$,
        \item [(c)] $\bullet_1$ $\mathcal M_\bft = \mathcal M_\gs$;

$\bullet_2$ if $a\in\mathcal M_\gs\setminus\clU_\gs$ then $M_{\bft,a} = \{\ell\}\times M_{\gs,a}$;

$\bullet_3$ if $a\in\clU_\gs$ then $M_{\bft,a} = M_{\gs,a}$,
        \item [(d)] $A_\bft = A_\bfu\times V_{(\ell,k)\times\gs}$, but if $\ell+1 = k$ let $A_\bft = A_\bfu$,
        \item [(e)] if $\nu\in{}^{(M_\bft)}\Lambda$ then $A^\bft_\nu = \{(\varrho,a) : \varrho\in V_{(\ell,k)\times\gs}$ and $a\in A^\bfu_{\varrho_*\cup\nu\cup\varrho}\}$, but if $\ell+1 = k$ then $A^\bft_\nu = \{a\in A_\bfu : a\in A^\bfu_{\varrho_*\cup\nu}\}$;
    \end{enumerate}
    \item [(C)] (a) $\bft$ indeed belongs to $\bfS$,

    (b) $\gs$ is isomorphic to $\gt = \mathrm{pre}(\bft)$, i.e.\ $\pi(b) = (\ell,b)$ for $b\in M_\gs$, $\hat\pi(\eta) = \nu$ for $\eta\in V_\gs$, where $\nu\in V_\gt$ and $(\forall b\in M_\gs)(\eta(b) = \nu(\pi(b)))$,

    (c) $\min(\bft)\ge\eps$,

    (d) if $Y\subseteq V_\gt$ and $\bigcap\{A^\bft_\nu : \nu\in Y\}\ne\emptyset$ then for some $\varrho\in V_{(\ell,k)\times\gs}$ we have $\bigcap\{A^\bfu_{\varrho_*\cup\nu\cup\varrho} : \nu\in Y\}\ne\emptyset$.
\end{enumerate}
\end{claim}

\begin{PROOF}{\ref{d65}}
Note that if $\LL\varrho_i : i<\ell\RR$ is defined we can choose $\bft_\ell$ and all is fine except possibly (C)(c). So by induction on $\ell<k$ we try to choose $\varrho_\ell$ such that
\begin{enumerate}
    \item [$(*)_\ell$] (a) $\varrho_\ell\in V_{\{\ell\}\times\gs}$,

    (b) let $\varrho_{<\ell} = \bigcup\{\varrho_h : h<\ell\}$, $\varrho_{\le\ell} = \varrho_{<\ell}\cup\varrho_\ell$,

    (c) $\eps_*+(\ell+1)\zeta\le\frac1{|V_{(\ell,k)\times\gs}|}\sum\big\{\mu_\bfu(A^\bfu_{\varrho_{\le\ell}\cup\eta}) : \eta\in V_{(\ell,k)\times\gs}\big\}$.
\end{enumerate}
We cannot succeed as the left side in $(*)_\ell$(c) is always $\le1$.
If we are stuck in $\ell$, i.e.\ we cannot choose $\varrho_\ell$, then easily $\bft_\ell$ is as required.
[Why? If we are stuck at $\ell$, some $\nu\in V_{\{\ell\}\times\gs}$ has $\mu(A^{\bft_\ell}_\nu)<\eps$. As the average over the $|V_\gs|$ values of $\nu$ is $\ge\eps_*+\ell\zeta$, the others average $\ge\eps_*+\ell\zeta+\frac{\eps_*-\eps}{|V_\gs|-1}$, i.e.\ we can choose $\varrho_\ell$ with increment $\zeta$.]
\end{PROOF}

\bn
\begin{observation}\label{d67}
1) The bound in \ref{d65} is one exponentiation, that is $k\le\big\lceil\frac1{\eps_*-\eps}\big\rceil\cdot|\Lambda|^{|M_\gs|}$ and $|M_\gu| = k\cdot|M_\gs|$.

2) The bound in \ref{d40} is $3+2r$ exponentiations, starting with $2r+|\Lambda|+\lceil\frac1\eps\rceil$.
\end{observation}

\bn
\begin{theorem}\label{d70}
$\HJ_{\eps_*}(\Lambda,r)$ is $\le m$ when
\begin{enumerate}
    \item [$\boxplus$] (a) $n = |\Lambda|$, $m = m_{|\Lambda|}$, $1\le t<r$, and e.g.\ $\eps = \eps_*/2$,

    (b) $\Lambda = \{\alpha,\beta_0,\dots,\beta_{n-2}\}$ with no repetitions so $|\Lambda| = n$, and $\Lambda_i = \{\alpha,\beta_j : j<i\}$,

    (c) $\bar\eps = \LL\eps_i : i<n-1\RR$, $\eps_i = \exp_\eps(\exp_3(i))$ for $i<|\Lambda|-1$, and $\eps_0 = \eps$,

    (d) $\LL m_i : i<n\RR$ is defined by induction by:
    \begin{itemize}
        \item [$\bullet_1$] $m_0 = r$,
        \item [$\bullet_2$] $(m_{i+1},m_i,\Lambda_i,\alpha,\beta_i,\eps_i,t)$ are like $(m,r,\Lambda,\alpha,\beta,\eps,t)$ in \ref{d40}, for $i<n-1$,
        \item [$\bullet_3$] $(m_{|\Lambda|},m_{|\Lambda|-1},\eps_*,\eps)$ are like $(m,r,\eps_*,\eps)$ in \ref{d65}.
    \end{itemize}
\end{enumerate}
\end{theorem}

\begin{PROOF}{\ref{d70}}
We need one exponentiation for applying \ref{d65}, e.g.\ choosing $\eps = \eps_*/2$, and we induct on $i<n-2$ dealing with $(\alpha,\beta_i)$ in the $i$-th step.

In details, let e.g.\ $\eps = \eps_*/2$ and let $\Lambda_i = \{\alpha,\beta_j : j<i\}$ for $i\le|\Lambda|-1$. We choose $\LL m_i : i\le|\Lambda|-1\RR$ by induction on $i$ such that $m_0 = r$ and $(m_{i+1},m_i,\Lambda_i,\alpha,\beta_i,\eps _i)$ $,t)$ are like $(|M_\gs|,r,\Lambda,\alpha,\beta,\eps ,t)$ in \ref{d40}.

Then choose $m_n = m$ such that $(m_{|\Lambda|},m_{|\Lambda|-1},\eps_*,\eps)$ are like $(m,r,\eps_*,\eps)$ in \ref{d65}. We choose
\begin{enumerate}
    \item [$(*)_{1}$] $\gs = (m,=)$, $\gs_* = \gs_{n-1}$.
\end{enumerate}
For $i\le|\Lambda|-1$ we choose
\begin{enumerate}
    \item [$(*)_{2}$] $\gs_i\in\gS_{\alpha,\beta_i,\Lambda}$ is defined for $i<n-1$ by

    (a) $|M_{\gs_i}| = m_{i+1}$,

    (b) $|\mathcal M_{\gs_i}| = m_i$,

    (c) $\clU_{\gs_i} = \emptyset$,

    (d) $\LL M_{\gs_i,a} : a\in\mathcal M_{\gs_i}\RR$ is as in \ref{d40}.
\end{enumerate}
We also define
\begin{enumerate}
    \item [$(*)_{3}$] $\gt_i\in\gS_{\alpha,\beta_i}$ is defined for $i<n-1$ by

    (a) $|M_{\gt_i}| = m_i = |\mathcal M_{\gs_i}|$,

    (b) $\clU_{\gt_i} = \mathcal M_{\gt_i}$ so $|M_{\gt_i,a}| = 1$ for $a\in\mathcal M_{\gt_i}$.
\end{enumerate}
To prove the claim we are given $X$:
\begin{enumerate}
    \item [$\oplus_1$] $X\subseteq{}^m\Lambda$, $|X|/|{}^m\Lambda|\ge\eps _*$.
    \item [$\oplus_2$] Let $\bfs\in\bfS$ be defined as $(m,=,A_\bfs,\bar A_\bfs)$ where $A_\bfs = \{0\}$, $A_{\bfs,\eta} = \{0\}$ if $\eta\in X$, $= \emptyset$ if $\eta\notin X$.
    \item [$\oplus_3$] Clearly $\av(\bfs)\ge\eps _*$.
    \item [$\oplus_4$] As $\gs\to_{\rm avm}(\gs_*)_{\eps_*,\eps}$ we can find $\bfs_*$ expanding $\gs_*$ embeddable into $\bfs$, that is isomorphic to some $\bfs'\in\Sub_\bfo(\bfs)$ such that $\min(\bfs_*)\ge\eps$.
    \item [$\oplus_5$] Find $\bfs_{n-1}$ expanding $\gs_{n-1}$ with $(A_{\bfs_{n-1}},\bar A_{\bfs_{n-1}}) = (A_{\bfs_*},\bar A_{\bfs_*})$, justified as $V_{\gs_{n-1}} = V_{\gs_*} = {}^{(m_{n-1})}\Lambda$.
    \item [$\oplus_6$] By downward induction on $i\le|\Lambda|-1$ we choose $\bfs_i$, $\bfo_i$, $\mathbf h_i$, $\bft_i$ such that

    (a) $\bfs_i$ expands $\gs_i$ and $\min(\bfs_i)\ge\eps _i$,

    (b) $\bft_i$ expands $\gt_i$ and $\min(\bft_i)\ge\eps _{i+1}$,

    (c) if $\bfs_i$ is defined (e.g.\ for $i = n-1$) then $\bft_i$ is isomorphic to some $\bfu_i\in\Sub _{\mathbf h_i,\bfo_i}(\bfs_i)$,

    (d) if $\bft_{i+1}$ is defined then $\bfs_i$ is expanding $\gs_i$ by:
    \begin{itemize}
        \item [($\alpha$)] $A_{\bfs_i} = A_{\bft_{i+1}} = A_{\bfs_{i+1}}$,
        \item [($\beta$)] if $\eta\in V_{\gs_i}$ then, recalling $M_{\gs_i} = M_{\gt_{i+1}}$ and $\Lambda_i = \Lambda_{i+1}\setminus\{\beta_i\}$:
        \item [$\circledast$] we have $\eta\in V_{\gs_i}$, $\beta_i\notin\mathrm{Rang}(\eta)$; we define $\bfo_\eta\in\mathrm{sob}(\gs_{i+1}) = \ob(\gt_{i+1})$ by: $\bullet$ $\mathcal N_{\bfo_\eta} = \eta^{-1}\{\alpha\}$, $\bullet$ $\rho_{\bfo_\eta} = \eta\rest(M_{\gt_i}\setminus\mathcal N_{\bfo_\eta})$, $\bullet$ $\cE_{\bfo_\eta}$ is the equality on $\mathcal N_{\bfo_\eta}$,
        \item [($\gamma$)] lastly $A^{\bfs_i}_\eta = \bigcap\{A^{\bfs_{i+1}}_\varrho : \varrho\in\Omega_{\bfo_\eta}\}$, recalling $\Omega_{\bfo_\eta}\subseteq V_{\gs_{i+1}}$.
    \end{itemize}

    Why can we carry the induction? $\bfs_{n-1}$ was already chosen. We can choose $\bft_i$ because $\gs_i\to_{\min}(\gt_i)_\eps$. Lastly having $\bft_{i+1}$, $\bfs_i$ has already been chosen then, and $\oplus_6$(d). That is, $\bfs_i$ is as required. The main point is ``$\bfs_i$ is $(\eps ,t)$-good'', which holds by our choices.]
    \item [$\oplus_7$] We can choose $\bfo_i\in\ob(\bfs_i)$ by induction on $i$ such that $\mu_{\bfs_i}(A^{\bfs_i}_{\Omega^+_{\bfo_i}})\ge\exp_\eps\exp_3(m_i)$, see \ref{d8}(3)$\bullet_4$ and \ref{d8}(6).

    [Why? For $i = 0$ this is trivial. For $i = j+1$, see the choices of $\bft_{j+1}$ and $\bfs_i$.]
\end{enumerate}
\end{PROOF}

\bn
\begin{conclusion}\label{d73}
The bound on $\HJ_\eps(\Lambda,r)$ is in $\clE^5$.
\end{conclusion}

\begin{PROOF}{\ref{d73}}
By \ref{d67} and \ref{d70}.
\end{PROOF}

\providecommand{\href}[2]{#2}
\section*{References}
\noindent [Cor69a] K.~Corr{\'a}di, \emph{{A remark on the number of common elements of sets (in Hungarian)}}, Mat. Lapok \textbf{20} (1969), 351--353.\par\smallskip
\noindent [Cor69b] Kereszt{\'e}ly Corr{\'a}di, \emph{{Problem at Schweitzer competition}}, Matematikai Lapok \textbf{20} (1969), 159--162.\par\smallskip
\noindent [FurKat91] H.~Furstenberg and Y.~Katznelson, \emph{A density version of the {H}ales-{J}ewett theorem}, J. Analyse Math. \textbf{57} (1991), 64--119.\par\smallskip
\noindent [GolMirShe:1182] Mohammad Golshani, Mostafa Mirabi, and Saharon Shelah, \emph{{[Sh:1182] From Shelah's block-content to Hales-Jewett}}, \href{https://arxiv.org/abs/1912.08643}{arXiv: 1912.08643}.\par\smallskip
\noindent [GraRotSpe90] Ronald Graham, Bruce~L. Rothschild, and Joel Spencer, \emph{Ramsey theory}, {Wiley -- Interscience Series in Discrete Mathematics}, Wiley, New York, 1990, 2nd edition.\par\smallskip
\noindent [HalJew63] A.W. Hales and R.I. Jewett, \emph{{Regularity and positional games}}, Trans. Amer. Math. Soc. \textbf{106} (1963), 222--229.\par\smallskip
\noindent [Juk11] S.~Jukna, \emph{{Extremal Combinatorics, with Applications in Computer Science, 2nd ed.}}, Springer, 2011.\par\smallskip
\noindent [Pol12] D.~H.~J. Polymath, \emph{A new proof of the density {H}ales-{J}ewett theorem}, Ann. of Math. (2) \textbf{175} (2012), no.~3, 1283--1327.\par\smallskip
\noindent [She:329] Saharon Shelah, \emph{{[Sh:329] Primitive recursive bounds for van der Waerden numbers}}, J. Amer. Math. Soc. \textbf{1} (1988), no.~3, 683--697, \href{https://www.ams.org/mathscinet-getitem?mr=929498}{MR929498} \href{https://doi.org/10.2307/1990952}{DOI: 10.2307/1990952}.\par\smallskip
\noindent [She:679] Saharon Shelah, \emph{{[Sh:679] A partition theorem}}, Sci. Math. Jpn. \textbf{56} (2002), no.~2, 413--438, \href{https://arxiv.org/abs/math/0003163}{arXiv: math/0003163} \href{https://www.ams.org/mathscinet-getitem?mr=1922806}{MR1922806}.\par\smallskip
\end{document}